\documentclass[10pt]{amsart}
\usepackage[foot]{amsaddr}

\usepackage{color}
\usepackage{multirow}
\usepackage{fullpage}

\usepackage{amsmath}
\usepackage{amssymb}
\usepackage{tikz}
\usepackage{tikz-cd}
\usepackage{diagbox}
\usepackage{lineno}

\newtheorem{proposition}{Proposition}
\newtheorem{corollary}{Corollary}
\newtheorem{remark}{Remark}

\newcommand{\algnd}[1]{\begin{aligned} #1 \end{aligned}}
\newcommand{\eqn}[1]{\begin{equation} #1 \end{equation}}
\newcommand{\eqns}[1]{\begin{equation*} #1 \end{equation*}}

\renewcommand{\div}{\operatorname{\nabla\cdot}}
\newcommand{\hdiv}{{\mathbf H}(\operatorname{div})}

\newcommand{\hdivz}{{\mathbf H}_0(\operatorname{div})}
\newcommand{\hhdivz}{{\mathbf H}_{0,h}(\operatorname{div})}

\newcommand{\norm}[1]{\left\Vert#1\right\Vert}

\newcommand{\inner}[2]{\left( #1, #2\right)}
\newcommand{\trinner}[2]{\left\langle #1, #2 \right\rangle}

\newcommand{\inv}[1]{#1^{-1}}
\newcommand{\grad}{\nabla}

\newcommand{\dual}[1]{#1^{\prime}}

\newcommand{\jump}[1]{\ensuremath{[\![#1]\!]} }
\newcommand{\avg}[1]{\ensuremath{\left\{\!\left\{#1\right\}\!\right\}} }

\newcommand{\reals}{\mathbb{R}}		

\newcommand{\uu}{\mathbf{u}}

\renewcommand{\SS}{\mathbf{S}}
\newcommand{\zz}{\mathbf{z}}
\newcommand{\ff}{\mathbf{f}}
\newcommand{\vv}{\mathbf{v}}
\newcommand{\VV}{\mathbf{V}}
\newcommand{\ww}{\mathbf{w}}
\newcommand{\ZZ}{\mathbf{Z}}

\newcommand{\LL}{\mathbf{L}}

\newcommand{\TT}[1]{\textcolor{black}{#1}}
\newcommand{\kent}[1]{\textcolor{black}{#1}}
\newcommand{\miro}[1]{\textcolor{black}{#1}}
\newcommand{\trygve}[1]{\textcolor{black}{#1}}

\newcommand{\Div}{\nabla\cdot}

\newcommand{\foralls}{\ensuremath{\forall\,}}

\def \lif {L_h}
\def \lif {\Theta_h}

\begin{document}


\title{An Observation on the Uniform 
Preconditioners for the mixed Darcy problem}

\author{Trygve B{\ae}rland}
\address[T.~B{\ae}rland]{Dept.~of Mathematics.~Univ.~of Oslo.~Oslo, Norway.}
\email[T.~B{\ae}rland]{trygveba@math.uio.no}

\author{Miroslav Kuchta}
\address[M.~Kuchta]{Dept.~of Num.~Anal.~and Sci.~Comput.~Simula Res. Lab. 
Fornebu, Norway.}
\email{miroslav@simula.no}

\author{Kent-Andre Mardal}
\address[K.-A.~Mardal]{
Dept.~of Mathematics.~Univ.~of Oslo.~Oslo, Norway and 
Dept.~of Num.~Anal.~and Sci.~Comput.~Simula Res.~Lab.~Fornebu, Norway.}
\email{kent-and@math.uio.no}

\author{Travis Thompson}
\address[T.~Thompson]{Dept.~of Num.~Anal.~and Sci.~Comput.~Simula Res.~Lab.~Fornebu, 
Norway. Present: Mathematical Institute, Oxford.}
\email{travis.thompson@maths.ox.ac.uk}

\date{\today}

\maketitle 
\begin{abstract}
When solving a multi-physics problem one often decomposes a monolithic system 
into simpler, frequently single-physics, subproblems.  A comprehensive solution 
strategy may commonly be attempted, then, 
by means of combining strategies devised for 
the constituent subproblems.  When decomposing the monolithic problem, however, 
it may be that requiring a particular scaling for one subproblem enforces an 
undesired scaling on another.  In this manuscript we consider the $\hdiv$-based mixed 
formulation of the Darcy problem as a single-physics subproblem; the hydraulic 
conductivity, $K$, is considered intrinsic and not subject to any rescaling.  
Preconditioners for such porous media flow problems in mixed form 
are frequently based on $\hdiv$ preconditioners rather than the pressure Schur 
complement. We show that when the hydraulic conductivity, $K$, is small the 
pressure Schur complement can also be utilised  for $\hdiv$-based 
preconditioners. The proposed 
approach employs an operator preconditioning framework to 
establish a robust, $K$-uniform block preconditioner.  The mapping properties 
of the continuous operator \trygve{is} a key component in applying the 
theoretical framework point of view\trygve{. As}  such, a main challenge addressed 
here is establishing a $K$-uniform inf-sup 
condition with respect to appropriately weighted Hilbert 
intersection\trygve{-} and sum spaces. 
\end{abstract}

\section{Introduction}
In this paper we will consider the mixed formulation of the Darcy problem 
of the form
\begin{eqnarray}
\label{mixeddarcy:1}
\frac{1}{K} \uu - \nabla p &=& \ff, \quad \mbox{ in } \Omega, \\ 
\label{mixeddarcy:2}
\nabla \cdot \uu &=& g, \quad \mbox{ in } \Omega,  
\end{eqnarray}
equipped with suitable boundary conditions. The variables $\uu$ and $p$ 
represent the fluid flux and pressure, respectively, and $K=\kappa/\mu_f$ 
denotes the hydraulic conductivity where $\kappa$ is the material permeability 
and $\mu_f$ the fluid viscosity. 

In this manuscript we pursue two purposes.  
The first purpose is to establish a uniform-in-$K$ inf-sup condition for 
\eqref{mixeddarcy:1}-\eqref{mixeddarcy:2}.  
Our second, and primary, purpose is the construction of some efficient 
block-diagonal preconditioners, for \eqref{mixeddarcy:1}-\eqref{mixeddarcy:2}, 
exhibiting robustness for $K\in (0,1)$.  These two objectives are connected. 
Indeed the operator preconditioning framework \cite{mardal2011preconditioning} 
\trygve{uses} the former stability result to provide for the latter \trygve{a} robust, 
block-diagonal preconditioner realization.  

In general, the framework approach is predicated on establishing a 
well-posedness result for the continuous problem 
in $K$-weighted Sobolev spaces.  For basic linear systems of the form  
\begin{equation}\label{eqn:primary-linear-system}
  \left[\begin{array}{cc} A & B^T \\ B & 0 \end{array} \right]
\left[\begin{array}{c} x  \\ y \end{array} \right] = 
\left[\begin{array}{c} c  \\ d \end{array} \right],   
\end{equation}
there are, generally speaking, two common approaches for constructing 
block-diagonal preconditioners. \trygve{One} may utilize a Schur complement for the 
first unknown or, alternatively, the second unknown.  That is, the structural 
options for the preconditioner are: 
\begin{equation}\label{eqn:structural-options}
\left[\begin{array}{cc} A^{-1} & 0 \\ 0 & (B A^{-1} B^T)^{-1} \end{array} 
\right]
\quad \mbox{and} \quad  
\left[\begin{array}{cc} (A+B^T B)^{-1} & 0 \\ 0 & X? \end{array} \right] . 
\end{equation}
The first approach results in three distinct, unit-sized eigenvalues~
\cite{murphy2000note} for any $0<K<\infty$.  On the other hand,  to the 
authors knowledge, only partial explanations have been offered for the second 
approach for the mixed Darcy problem.  In \cite{arnold1997preconditioning}, 
$\hdiv$ preconditioners were constructed and applied to 
\eqref{mixeddarcy:1}-\eqref{mixeddarcy:2} in the case $K=1$ for which $X$ 
coincides with the inverse of a potentially diagonalized mass matrix.  The 
same authors also developed multilevel methods for weighted $\hdiv$ spaces in 
\cite{arnold2000multigrid}, but did not discuss the corresponding appropriate 
scaling for a Darcy problem.  In \cite{powell2003optimal} the consequences of 
$K$-scaling, also of spatially varying $K$, were studied for both approaches; 
it was shown that the eigenvalues of the preconditioned system were affected 
significantly by $K$ and good results were obtained only when a proper 
rescaling was used. 

If we now consider a scaled version \eqref{eqn:primary-linear-system} we 
arrive at a system with the general form %
\begin{equation}\label{eqn:primary-scaled-linear-system}
\left[\begin{array}{cc} \alpha A & B^T \\ B & 0 \end{array} \right]
\left[\begin{array}{c} x  \\ y \end{array} \right] = 
\left[\begin{array}{c} c  \\ d \end{array} \right]. 
\end{equation}
In \cite{hong2017parameter,vassilevski2013block,vassilevski2014mixed} it was 
suggested that problems of the form of \eqref{eqn:primary-scaled-linear-system} 
may be preconditioned efficiently by a scaling of the second structural option 
of \eqref{eqn:structural-options}; that is, preconditioners structurally based 
on
\begin{equation}
\label{their_prec:a}
\mathcal{B}_1 = \left[\begin{array}{cc} (\alpha A + \alpha B^T B)^{-1}  & 0 \\ 
0 & \TT{\left(\frac{1}{\alpha} I\right)^{-1} } \end{array} \right].
\end{equation}
This view was advanced in the aforementioned work, to develop preconditioners 
for the Biot and Brinkman problems, where the Darcy part was preconditioned 
by the approach of \eqref{their_prec:a}.  The Brezzi conditions for the Darcy 
problem in the weighted norms corresponding to the choice of 
\eqref{their_prec:a}, that is $\frac{1}{\sqrt{K}}\hdivz\times\sqrt{K}L_0^2$, 
follow directly from the unscaled version by a simple scaling.  Furthermore, 
the approach is closely related to the augmented Lagrangian approach which 
was investigated in \cite{benzi2006augmented,greif2006preconditioners} 
and used successfully for the Oseen and Maxwell type problems, respectively.  

Given the form of \eqref{their_prec:a} we are motivated to ask the following 
question: is the $\alpha$-scaling of $B^T B$, in the top-left block, necessary?  
The motivation for this question comes from considering 
\eqref{mixeddarcy:1}-\eqref{mixeddarcy:2} as part of a multi-physics problem.  
For instance, as a single-physics problem, the isolated case of 
$K\rightarrow 0$, in \eqref{mixeddarcy:1}-\eqref{mixeddarcy:2}, is not 
necessarily intrinsically interesting---a simple scaling resolves the issue of 
a vanishing $K$.  However, in a multi-physics setting, solution algorithms
are often constructed via decomposition into single-physics subproblems. Thus, 
requiring a certain scaling of one of the single-physics sub-problems, as we 
see here with the Darcy problem, may enforce an undesired scaling on other 
single-physics problems within the multi-physics system. 

To investigate this question we consider an alternative to 
\eqref{their_prec:a}; namely   
\begin{equation}
\label{our_prec}
\mathcal{B}_2 = \left[\begin{array}{cc} (\alpha A + B^T B)^{-1}  & 0 \\ 
0 & I\TT{^{-1}} +  (B(\alpha A)^{-1} B^T)^{-1}    
\end{array} \right] . 
\end{equation}

In the current work we construct block preconditioners based on 
the operator preconditioning framework \cite{arnold1997preconditioning,
mardal2011preconditioning}.  In the continuous case, for the Darcy problem,  
the preconditioner takes the particular form, of \eqref{our_prec}, given by 
\[
\mathcal{B} = \left[ \begin{array}{cc} (K^{-1} - \nabla \div)^{-1} & 0 \\ 
           0 & \TT{I^{-1}} + (-\div(K\grad))^{-1} \end{array} \right]. 
\]
It will be shown that preconditioners based on $\mathcal{B}$ are robust with 
respect to both the mesh size and the permeability $K$.  We will also consider 
permeabilites with jumps and anisotropy in our numerical experiments where we 
compare our proposed approach, \eqref{our_prec}, to the previous approach of 
\eqref{their_prec:a}. The remainder of this manuscript is organized as follows: 
Section \ref{sec:prelim} introduces the necessary notation and basic results; 
Section \ref{sec:uniform_conts} discusses the continuous uniform stability 
(inf-sup) condition and the resulting continuous preconditioner; Section 
\ref{sec:uniform_disc} addresses the corresponding discrete case of each; in 
Section \ref{sec:numerical-experiments} \trygve{ the preconditioners proposed in this paper are validated through }numerical experiments. \trygve{We also use this section to explore some cases of practical interest which are not covered directly by the theory of Sections \ref{sec:uniform_conts} and \ref{sec:uniform_disc}.} 
Finally, Section \ref{sec:conclusion} offers some concluding remarks. 
\trygve{In closing,  we mention} that the theory of Sections 
\ref{sec:uniform_conts}\trygve{--}\ref{sec:uniform_disc} assume $K \in (0,1)$ is an 
arbitrary, but fixed, constant. 

\section{Preliminaries}\label{sec:prelim}
Let $\Omega$ be a bounded, connected Lipschitz domain in $\reals^n,\, n=2,3$.
Then $L^2=L^2(\Omega)$ denotes the 
space of square integrable functions, whereas
$H^k=H^k(\Omega)$ denotes the Sobolev space of functions with all derivatives up to order $k$ in $L^2$.
The spaces $L^2/\reals$ and $H^1/\reals$  contain functions in $L^2$ and $H^1$, respectively, with zero mean value. 
Vector valued functions, and Sobolev spaces of vector valued functions, are denoted
by boldface. The space $\hdiv$ contains functions in  ${\bf L}^2$ with 
divergence in $L^2$ and the subspace of functions $\uu \in \hdivz$ are those 
with zero normal trace. The notation $(\cdot, \cdot)$ is used for the $L^2$ inner 
product and analogously for vector fields. The norm corresponding to the $L^2$ inner 
product is expressed with the canonical double-bar $\norm{\cdot}_{L^2}$.  
The duality pairing between a real Hilbert space $X$ and its dual 
$\dual{X}$ is $\trinner{\cdot}{\cdot}$.  
 \TT{Suppose that $X$ is a Hilbert space and that $\alpha > 0$ 
is a fixed real value.  Then we denote by $\alpha X$ the Hilbert space whose 
elements coincide with the elements of $X$ and with norm %
\[
\norm{f}_{\alpha X}^2 = \alpha^2 \norm{f}_X^2.
\]
}


For two 
\TT{Hilbert} spaces, $X$ and $Y$, we denote by $\mathcal{L}(X,Y)$ the
space of bounded linear maps from $X$ to $Y$. 
In the subsequent analysis we employ both the intersection and sum of two 
\TT{Hilbert} spaces $X$ and $Y$. These composite spaces are formally defined as 
follows: Let $X$ and $Y$ be two 
\TT{Hilbert} spaces \TT{which are both subspaces of some} 
larger 
\TT{Hilbert} space\trygve{. The} intersection and sum space are \TT{defined, respectively, as}:
\TT{%
\[
 X\cap Y = \left\{ f \, \vert \, f\in X, \text{ and } f \in Y \right\}, %
\quad X+Y = \left\{ f + g \, \vert \, f\in X, \text{ and } g \in Y \right\}.
\]
In this manuscript we will be concerned with the case where $X$ and $Y$ are 
Hilbert spaces. In this case $X\cap Y$ and $X+Y$ are also Hilbert 
spaces with respect to the norms:
\begin{equation}\label{eqn:capnorm}
	\|z\|^2_{X \cap Y} = \|z\|_X^2 + \|z\|_Y^2   
\end{equation}}
and \TT{%
\begin{equation}\label{eqn:sumnorm}
	\|z\|^2_{X + Y} = \inf_{\substack{z=x+y\\x\in X, y\in Y}}\|x\|_X^2 + \|y\|_Y^2.  
\end{equation}
}


\noindent In addition \TT{\cite[Theorem 2.7.1]{lofstrom1976}} if $X \cap Y$ is 
dense in both $X$ and $Y$, then
\eqn{
	\label{eq:sum-space-dual}
	\dual{\left( X + Y \right)} = \dual{X} \cap \dual{Y}.
}
\TT{
\begin{remark}
In \cite{lofstrom1976} $X$, $Y$, $X\cap Y$ 
and $X+Y$ are Banach spaces. Again, see \cite[Lemma 2.3.1]{lofstrom1976}.  The norms 
of $X\cap Y$ and $X+Y$ are given 
explicitly by the alternative definitions 
\begin{equation}\label{eqn:bannorms:lof}
	\|z\|_{X \cap Y} = \max \left(\|z\|_X, \|z\|_Y\right), \text{ and } 
\quad \|z\|_{X + Y} = \inf_{\substack{z=x+y\\x\in X, y\in Y}}\|x\|_X + \|y\|_Y,  %
\end{equation}
which is different from the norms above \eqref{eqn:capnorm} and \eqref{eqn:sumnorm}.  
In our context, where we aim to derive preconditioners, the norms  
\eqref{eqn:capnorm} and \eqref{eqn:sumnorm} are, however, more convenient.
We therefore detail the equivalence between the two different definitions. 
In particular, 
\[
\frac{1}{2}\left(\norm{z}_X^2 + \norm{z}_Y^2 \right) \leq %
\max \left(\|z\|_X, \|z\|_Y\right) ^2 \leq \norm{z}_X^2 + \norm{z}_Y^2,
\]
and further 
\[
\inf_{\substack{z=x+y\\x\in X, y\in Y}}\|x\|_X^2 + \|y\|_Y^2 %
\leq \left(\inf_{\substack{z=x+y\\x\in X, y\in Y}}\|x\|_X + \|y\|_Y\right)^2 %
\leq 2 \inf_{\substack{z=x+y\\x\in X, y\in Y}}\|x\|_X^2 + \|y\|_Y^2
\]
Hence, the norms employed here, \eqref{eqn:capnorm} and \eqref{eqn:sumnorm}, are equivalent with the norms in~\cite{lofstrom1976}.
Finally, we remark that in the more general 
Banach space setting \eqref{eq:sum-space-dual} is understood as an identification 
under an isometry, c.f. \cite[Theorem 2.7.1]{lofstrom1976}, with respect to 
the norms of \eqref{eqn:bannorms:lof}.  In the Hilbert case the map providing 
the identification \eqref{eq:sum-space-dual} is, instead, an isomorphism with 
respect to the (equivalent) norms \eqref{eqn:capnorm} and \eqref{eqn:sumnorm}.  
\end{remark}
}
Now suppose that 
$\{X_1, X_2 \}$ and $\{Y_1, Y_2 \}$ are pairs of \TT{Hilbert} spaces
such that both elements of each pairing are \TT{subspaces of} \trygve{some} \TT{larger Hilbert space}.
If $T$ is a bounded linear operator from $X_i$ to $Y_i$ for $i = 1,\, 2$, then
\eqn{
	\label{eq:operator-sum-spaces}
	T \in \mathcal{L}(X_1 \cap X_2, Y_1\cap Y_2) \cap \mathcal{L}(X_1 + X_2, Y_1 + Y_2),
}
and in particular
\eqns{
	\label{eq:sum-space-operator-norm}
	\norm{T}_{\mathcal{L}(X_1 + X_2, Y_1 + Y_2)} \leq \max \left(\norm{T}_{\mathcal{L}(X_1,Y_1)}, \norm{T}_{\mathcal{L}(X_2,Y_2)} \right).
}

\section{Continuous stability and preconditioning}\label{sec:uniform_conts}

The weak formulation of \eqref{mixeddarcy:1}--\eqref{mixeddarcy:2} 
reads: Find $\uu \in \VV, p \in Q$ such that  
\begin{eqnarray}
a (\uu, \vv) + b (\vv, p) &=& (\ff, \vv), \quad \forall \vv \in \VV, \label{eqn:mixeddarcy:1:w}\\   
b(\uu, q) &=& (g, q), \quad \forall q \in Q ,\label{eqn:mixeddarcy:2:w} 
\end{eqnarray}
where  
\[
a(\uu, \vv) = \inner{K^{-1} \uu}{ \vv } \mbox{ and }   b(\uu, q) = \inner{\div \uu}{ q }.  
\]
We recall that $K\in(0,1)$ is considered an arbitrary but fixed constant. The corresponding coefficient matrix reads,   
\eqn{
	\label{eq:mixed-darcy-coefficient-matrix}
	\mathcal{A}\begin{bmatrix}
		\uu \\
		p 
	\end{bmatrix}
	:=
	\begin{bmatrix}
		K^{-1} & -\grad \\
		\div & 0 
	\end{bmatrix}
	\begin{bmatrix}
		\uu \\
		p 
	\end{bmatrix}
	=
	\begin{bmatrix}
		\ff \\
		g
	\end{bmatrix}.
}
\TT{%
The primary result of this section can now be stated: %
\begin{proposition}[A uniform-in-$K$ continuous stability result]\label{prop:uniform-inf-sup}
Consider the problem \eqref{eqn:mixeddarcy:1:w}-%
\eqref{eqn:mixeddarcy:2:w} where $\uu$ and $p$ are chosen, respectively, 
in the spaces: 
\begin{equation}\label{eqn:spaces}
\VV = K^{-1/2} \mathbf{L^2} \cap \hdivz, \quad Q = L^2/\reals + K^{1/2} (H^1 / \reals).
\end{equation}
Then the Brezzi conditions are satisfied uniformly in $K$. 
\end{proposition}
}


\begin{proof}
\TT{We begin with two comments: first, for applications to preconditioning, we 
are particularly interested in the existence of an inf-sup condition for 
$b(\vv,q)$ that is independent of the choice of $K\in(0,1)$.  We will prove 
that such a condition holds for the pair of spaces $\VV$ and $Q$. Second, we 
remark that the arguments below are similar to~\cite{mardal2013uniformly}, 
although the Sobolev spaces and bilinear forms involved are different.} We 
now verify all of the requisite Brezzi conditions for \eqref{eqn:mixeddarcy:1:w}--\eqref{eqn:mixeddarcy:2:w}. Towards 
this end let 
\[
\ZZ = \{ \uu \in \VV \ | \ b(\uu, q) = 0, \quad \forall q \in Q \} .   
\] 
Then clearly 
\[
a(\zz,\zz) = \|\zz\|^2_{K^{-1/2} \mathbf{L}^2}=\|\zz\|^2_{\VV}, \quad \forall \zz\in \ZZ  
\]
and hence coercivity of $a(\cdot, \cdot)$ over $\ZZ$ is established. Furthermore, the boundedness of 
$a(\cdot, \cdot)$ follows from 
\[
a(\uu,\vv) = (\uu, \vv)_{K^{-1/2}} \le \|\uu\|_{\VV} \|\vv\|_{\VV} \quad \forall \uu, \vv \in \VV.  
\]
The boundedness of $b$ follows from a decomposition argument. Let $q=q_0 + q_1$, $q_0\in L^2/\reals$ and $q_1 \in H^1 / \reals$ then \begin{eqnarray*}
b(\vv,q) &=& \inner{\nabla\cdot \vv}{q_0} - \inner{\vv}{\grad q_1} \\
       &=& \inner{\nabla\cdot \vv}{q_0} - \inner{K^{-1/2}\vv}{K^{1/2}\grad q_1} \\
       &\leq& \norm{\nabla\cdot \vv}_{L^2}\norm{q_0}_{L^2} + \norm{\vv}_{K^{-1/2}\LL^2}\norm{\grad q_1}_{K^{1/2}\LL^2} \\
       &\leq& \left(\inner{\inv{K}\vv}{\vv} + \norm{\nabla\cdot \vv}^2_{L^2} \right)^{\frac{1}{2}}\left(\norm{q_0}^2_{L^2} + \inner{K \grad q_1}{\grad q_1} \right)^{\frac{1}{2}}.
\end{eqnarray*}
Taking the infimum over all decompositions of $q$ yields the desired bound.   

To establish the uniform inf-sup condition we will demonstrate the existence of 
a linear operator $\SS$ which satisfies the following properties:
\begin{align}
&\SS \in \mathcal{L}(\dual{Q}, \VV),\label{eqn:S-properties-a}\\
&\norm{\SS}_{\mathcal{L}(\dual{Q}, \VV)}\text{ is independent of } K,\label{eqn:S-properties-b}\\
& \TT{\inner{\div \SS g}{\psi} = \trinner{g}{\psi}, \text{ for any }g \in \dual{Q} \text{ and } \psi \in Q.}\label{eqn:S-properties-c}
\end{align}
Suppose first that such an operator $\SS$ exists satisfying conditions 
\eqref{eqn:S-properties-a}--\eqref{eqn:S-properties-c}\trygve{, and} define the constant 
$C_{\SS} =  \norm{\SS}_{\mathcal{L}(\dual{Q}, \VV)}$.  Under these assumptions 
the inf-sup condition follows directly. To see this, take any $p \in Q$; 
it follows that:
\begin{equation}\label{eqn:uniform-inf-sup-qp-v}
	\sup_{\vv \in \VV} \frac{\inner{\div \vv}{p}}{\norm{\vv}_\VV} 
	\geq \sup_{g \in \dual{Q}} \frac{\inner{\div \SS g}{p}}{\norm{\SS g}_\VV} 
	= \sup_{g \in \dual{Q}} \frac{\trinner{g}{p}}{\norm{\SS g}_\VV} 
	\geq \inv{C_S} \sup_{g \in \dual{Q}}\frac{\trinner{g}{p}}{\norm{g}_{\dual{Q}}} 
	=\inv{C_S} \norm{p}_Q.
\end{equation}
We now 
show that such an operator $\SS$ can be realized as the solution of a suitable 
Poisson problem.  Moreover, the $K$-independence of the operator norm will 
follow from a suitable scaling between the primal and mixed formulations of the 
problem.  %
Towards this end let $g \in \dual{(H^1 / \reals)}$ 
and define $\phi \in H^1 / \reals$ to be the 
solution of the homogeneous Neumann problem, i.e., 
\eqn{
	\label{eq:S-def1}
	\inner{\nabla \phi}{\nabla \psi} = -\trinner{g}{\psi}, \quad \forall \psi \in H^1 / \reals.
}
Now define $\SS$ by $\SS g = \nabla \phi$; then \eqref{eq:S-def1} implies 
\eqn{
	\label{eq:S-mapping1}
	\SS \in \mathcal{L}(\dual{(H^1 / \reals)}, \LL^2).
}
When interpreted in the weak sense, \eqref{eq:S-def1} provides the 
identity \eqref{eqn:S-properties-c}; that is $\div \SS g = g$ where 
$\div:\LL^2\rightarrow \dual{(H^1/\mathbb{R})}$.  

For a second perspective \trygve{on} establishing the operator $\SS$, we now assume that 
$g \in \dual{(L^2/\reals)}$.  Consider then a mixed formulation of 
\eqref{eqn:mixeddarcy:1:w}--\eqref{eqn:mixeddarcy:2:w},   
with $K = 1$, given by: find $(\mathbf{r}, \phi) \in \hdivz \times L^2/\reals $ such that
\eqn{
	\label{eq:S-def2}
	\algnd{
		\inner{\mathbf{r}}{\mathbf{s}} + \inner{\div \mathbf{s}}{\phi} &= 0, \quad \forall \mathbf{s} \in \hdivz \\
		\inner{\div \mathbf{r}}{\psi} &= \trinner{g}{\psi}, \quad \forall \psi \in L^2/\reals.
	}
}
The above problem is classical\trygve{, and} as a result of its well-posedness we 
define $\SS g = \mathbf{r}$.  It is a straightforward exercise to establish that, 
for $g  \in L^2/\reals$, the above definition of $\SS g$ coincides with that of \eqref{eq:S-def1}.
Moreover we have that $\div \SS g = g$, from the second equation of \eqref{eq:S-def2},
and 
\eqn{
	\label{eq:S-mapping2}
	\SS \in \mathcal{L}(L^2/\reals, \hdivz).
}

The primary observation, then, is that a scaling argument now reveals that 
$\SS\in\mathcal{L}(Q',\VV)$ with operator norm independent of $K$.  To see this,
first scale \eqref{eq:S-mapping1} by $K^{-1/2}$; i.e.~ \eqref{eq:S-mapping1} implies
\[
	\SS\in\mathcal{L}(K^{-1/2}\dual{(H^1/\mathbb{R})}\,,\,K^{-1/2}\LL^2).
\]
Combine the above with \eqref{eq:S-mapping2} and use \eqref{eq:operator-sum-spaces} 
to get
\[
\SS\in\mathcal{L}\left(K^{-1/2}\dual{(H^1/\mathbb{R})}\cap L^2/\reals\,,\, K^{-1/2}\LL^2\cap\hdivz\right)=\mathcal{L}(\dual{Q},\VV),
\]
for any $K > 0$. 

In particular note that the operator norm of $\SS$ is independent of $K$. This
follows since the operator norms related to \eqref{eq:S-mapping1}, both before and 
after rescaling, and \eqref{eq:S-mapping2} are independent of $K$. 
We have shown that $\SS$, defined by \eqref{eq:S-mapping1}--\eqref{eq:S-mapping2}, 
satisfies the properties \eqref{eqn:S-properties-a}--\eqref{eqn:S-properties-c}.  
Thus, the desired uniform inf-sup condition is established, pursuant to the discussion
preceding \eqref{eqn:uniform-inf-sup-qp-v}. 
\end{proof}

\begin{remark}
Note that in the proof of Proposition \ref{prop:uniform-inf-sup}, above, we 
have made use of \eqref{eq:sum-space-dual}, and reflexivity, to establish 
that
\[
\dual{Q} = \dual{\left(L^2/\reals + K^{1/2}\left(H^1/\mathbb{R}\right)\right)} = L^2/\reals \cap K^{-1/2}\dual{\left(H^1/\mathbb{R}\right)}.
\]
\end{remark}

With a $K$-uniform stability established, the framework put forth in~\cite{mardal2011preconditioning}
details the construction of a robust and efficient preconditioner for $\mathcal{A}$ can be constructed
from preconditioners for the operators realizing the
$\VV$- and $Q$\trygve{-norms}.  The following proposition states precisely this: 

\TT{%
\begin{corollary}[A $K$-robust preconditioner for the continuous mixed Darcy problem]\label{prop:conts-precon}
The preconditioner defined by %
 \eqn{
	\label{eq:continuous-preconditioner}
	\mathcal{B} =
	\begin{bmatrix}
		(K^{-1}I - \grad \div)^{-1} & 0 \\
		0 & \TT{I^{-1}} + (-K\Delta)^{-1}
	\end{bmatrix}
}
provides a robust preconditioner for \eqref{eqn:mixeddarcy:1:w}--\eqref{eqn:mixeddarcy:2:w}
in the sense that the condition number of $\mathcal{B}\mathcal{A}$ is bounded uniformly in $K$. 
Here, $\mathcal{A}$ is given by \eqref{eq:mixed-darcy-coefficient-matrix}.
\end{corollary}
}

\begin{proof}
Preconditioners for the $\VV$ inner product are well-known, cf. e.g.~\cite{arnold2000multigrid,dobrev2019,hiptmair2007nodal,vassilevski2012pAMG,leeCS2016}.
\TT{Following the framework of \cite{mardal2011preconditioning} it therefore 
suffices to construct an operator realizing the $Q$-norm. The 
result can then be combined with the operator realizing the $\VV$-norm from 
the sources cited above}.
For the $Q$-norm, we begin by \TT{remarking that \eqref{eqn:sumnorm} is equivalent to the expression:}
\begin{equation}\label{eqn:Qnorm:conts}
	\norm{q}_Q^2 = \inf_{\phi \in H^1 /\reals} \left( \norm{q - \phi}_{L^2}^2 + K\norm{\grad \phi}_{L^2}^2 \right).
\end{equation}
\TT{In fact, the infimum on the right hand side is attained by the unique $\phi \in H^1/\reals$ solving the variational} problem:
\begin{equation}\label{eqn:qnorm:varform}
	\inner{\phi}{\psi} + K \inner{\grad \phi}{\grad \psi} = \inner{q}{\psi}, \quad \forall \psi \in H^1 / \reals.
\end{equation}
\TT{The above corresponds to the statement that $\phi$ is the solution of the 
elliptic problem %
\[ %
(I-K\Delta)\phi = q, \quad \grad \phi \cdot \mathbf{n}=0 \text{ on } \partial \Omega%
\] %
where $\mathbf{n}$ is the outward pointing normal to $\partial \Omega$, and 
recall that $q \in L^2/ \reals$.
We will make use of the expressions } $q - \phi = -K\Delta \phi$ 
\TT{and, equivalently,} $\phi = (I - K\Delta)^{-1}q$. 
\TT{Since $\phi$ 
minimizes \eqref{eqn:Qnorm:conts}, we have from the observations above that  %
\begin{align}%
\norm{q}_Q^2 &=  \norm{q - \phi}^2 + K\norm{\grad \phi}^2 = \inner{-K\Delta \phi}{-K\Delta \phi} + K\inner{\grad \phi}{\grad \phi} \\
& \quad=  \inner{-K\Delta \phi}{-K\Delta \phi} + \inner{-K\Delta \phi}{\phi} = \inner{(-K\Delta)\phi}{(I-K\Delta)\phi}\nonumber \\
& \quad= \inner{(-K\Delta)(I-K\Delta)^{-1}q}{q}. \nonumber
\end{align}%
}
%
The canonical preconditioner, $\mathcal{B}_Q: \dual{Q} \to Q$, corresponding 
to the above is
\eqns{
	\mathcal{B}_Q = \left[(-K\Delta) (I-K\Delta)^{-1}\right]^{-1} = I\TT{^{-1}} + (-K\Delta)^{-1}.
}
\TT{
According to  \cite{mardal2011preconditioning}\trygve{, a continuous and robust-in-$K$ preconditioner for}
\eqref{eq:mixed-darcy-coefficient-matrix}, is then 
\[
\mathcal{B} =
	\begin{bmatrix}
		\mathcal{B}_{\VV} & 0 \\
		0 & \mathcal{B}_Q
	\end{bmatrix}.
\] %
This is precisely \eqref{eq:continuous-preconditioner}.
}
\end{proof}

\section{Discrete Stability and preconditioning}\label{sec:uniform_disc}
In this section we describe the construction of a preconditioner for discretizations
based on Brezzi-Douglas-Marini (BDM) and Raviart-Thomas (RT) elements ~\cite{brezzi1985two,raviart1977mixed}.
Again we assume $K\in(0,1)$ is an arbitrary but fixed constant.  The discrete approach reflects many aspects of 
the continuous setting of section \ref{sec:uniform_conts}. 
However, due to the discontinuous polynomial nature of the pressure elements, we first define a 
discrete $H^1$-norm to establish the $Q$-norm in the discrete case.

Let $\mathcal{T}_h$  be a shape regular simplicial mesh defined on the bounded, 
Lipschitz domain $\Omega$ and let $r\geq 0$. Let $\VV_h$ be the 
$\hdiv$-conforming discrete space given by either the RT elements of order $r$ 
or the BDM elements of order $r+1$.  Define $Q_h$ to be the usual corresponding 
space of discontinuous, piecewise polynomials of order $r$.  Consider the 
discrete mixed Darcy problem given by: find $\uu_h \in \VV_h$ and $p_h \in Q_h$ 
such that
\begin{eqnarray}
a (\uu_h, \vv) + b (\vv, p_h) &=& (\ff, \vv), \quad \forall \vv \in \VV_h, \label{eqn:mixeddarcy:1:d}\\   
b(\uu_h, q) &=& (g, q), \quad \forall q \in Q_h.\label{eqn:mixeddarcy:2:d} 
\end{eqnarray}

The discrete $Q_h$-norm will be defined
in terms of a discrete gradient which is the negative $L^2$-adjoint of the $\div$ operator on $\VV_h$.
First, $\grad_h:Q_h\to\VV_h$ is defined by

 \begin{equation}\label{eqn:gradh}
	\inner{\grad_h q}{\vv} = - \inner{q}{\div \vv}. 
 \end{equation}
It is well-known, \cite{brezzi_fortin,Gatica_2014}, that with these particular choices of $Q_h$ and $\VV_h$, 
there is an $h$-independent constant $\beta > 0$ such that
\eqn{
	\label{eq:discrete-inf-sup}
	\sup_{\vv \in \VV_h}\frac{\inner{q}{\div \vv}}{\norm{\vv}_{\hdiv}} \geq \beta \norm{q}_{\LL^2}
}
for every $q \in Q_h$.
It follows that $\grad_h$ is injective and we can define the discrete $H^1$-norm on $Q_h$ via
\eqns{
	\norm{q}_{1,h} = \norm{\grad_h q}_{\mathbf{L^2}}.
} 
We denote the space $H_h^1$ as the set $Q_h$ equipped with the norm $\norm{\cdot}_{1,h}$
and the space $L^2_h$ as the set $Q_h$ equipped with the usual $L^2$-norm.  
The discrete analogue of the $Q$-norm, \TT{i.e.~the discrete form  of %
\eqref{eqn:Qnorm:conts} 
which is itself equivalent to 
the discrete form of \eqref{eqn:sumnorm}}
is given as
\eqn{
	\label{eq:Q_h-norm}
	\norm{q}_{Q_h}^2 = \inf_{\phi \in Q_h} \left(\norm{q - \phi}^2_{L^2} + K\norm{\phi}_{1,h}^2 \right).
}
\TT{The spaces $\hhdivz$ and $\mathbf{L}^2_{h}$ are defined analogously to $L^2_h$ above.  
The primary result of this section, which we now state, is the discrete 
analogue of Proposition \ref{prop:uniform-inf-sup}}.

\TT{%
\begin{proposition}[A uniform-in-$K$ discrete stability result]\label{prop:uniform-inf-sup}
Consider the problem \eqref{eqn:mixeddarcy:1:d}-%
\eqref{eqn:mixeddarcy:2:d} where $\uu_h$ and $p_h$ are chosen, respectively, 
in the spaces: 
\begin{equation}\label{eqn:spaces:disc}
\VV_h = K^{-1/2}\mathbf{L}^2_{h} \cap \hhdivz \quad \text{ and } \quad Q_h = L^2_h + K^{1/2} H^1_h,
\end{equation}
Then, under the corresponding norms, the Brezzi conditions are satisfied 
uniformly in $K$. 
\end{proposition}
}


\begin{proof}
Boundedness and coercivity of $a(\cdot,\cdot)$ and 
boundedness of $b(\cdot,\cdot)$ follows from the same arguments
put forth in section \ref{sec:uniform_conts}.    
Verifying a $K$-independent inf-sup condition will therefore conclude
the argument.  To accomplish this, a left-inverse of $\grad_h$ will be
constructed, satisfying appropriate bounds, allowing for a
similar argument to that of section \ref{sec:uniform_conts} 
for the operator $\SS$.  Let $\ZZ_h$ denote the discrete kernel 
of the $\div$ operator; i.e.~the set of $\vv_h \in \VV_h$ for which
\eqn{
\inner{\div\vv_h}{q_h} = 0,\quad\foralls q_h\in Q_h.
} 
From \eqref{eq:discrete-inf-sup} it follows (cf. \cite{Gatica_2014}) that $\grad_h:Q_h\to \ZZ_h^{\perp}$ is a linear bijection.  Furthermore, every $\vv\in\VV_h$ can be uniquely decomposed as 
\eqn{
	\label{eq:discrete-helmholtz}
	\vv = \grad_h r + \tilde{\vv},
}
where $r \in Q_h$, with $\grad_h r \in \ZZ_h^{\perp}$, and $\tilde{\vv} \in \ZZ_h$.  Since the spaces considered
for $\VV_h$ satisfy the relation $\div\VV_h\subset Q_h$ it follows that $\div\tilde{\vv}=0$, 
for every $\tilde{\vv}\in \ZZ_h$, and the decomposition \eqref{eq:discrete-helmholtz} is orthogonal with 
respect to both the $\hdiv$- and $\LL^2$ inner products. 

We now define the lifting operator $\lif: \VV_h \to Q_h$ by $\lif \vv = r$,
according to \eqref{eq:discrete-helmholtz}.
It is evident that $\lif \grad_h$ is the identity operator on $Q_h$ and that 
$\lif \tilde{\vv}=0$ for all $\tilde{\vv}\in \ZZ_h$.  Moreover, the inf-sup 
condition \eqref{eq:discrete-inf-sup} and the $\hdiv$-orthogonality of 
\eqref{eq:discrete-helmholtz} implies that
\eqns{
	\algnd{
		\norm{\lif \vv} &\leq \beta^{-1}\sup_{\ww \in \VV_h}\frac{\inner{\lif \vv}{\div \ww}}{\norm{\ww}_{\hdiv}} \\
		&= \beta^{-1}\sup_{\ww \in \VV_h}\frac{\inner{\grad_h \lif \vv}{\ww}}{\norm{\ww}_{\hdiv}} \\
		&\leq \beta^{-1}\sup_{\ww \in \VV_h}\frac{\inner{\vv}{\ww}}{\norm{\ww}_{\hdiv}} \\
		&= \beta^{-1}\norm{\vv}_{\dual{\hhdivz}},
	}
}
which means that $\lif \in \mathcal{L}(\dual{\hhdivz}, L^2_h)$.
From the $L^2$-orthogonality of \eqref{eq:discrete-helmholtz} we also have that
$\norm{\lif \vv}_{1,h} \leq \norm{\vv}$, which implies that $\lif \in \mathcal{L}(\mathbf{L}_h^2, H_h^1)$.
From these bounds on $\lif$, together with \eqref{eq:operator-sum-spaces}, we deduce that
\eqn{
	\label{eq:lifting-mapping-property}
	\lif \in \mathcal{L}\left(K^{1/2}\mathbf{L}_h^2 + \dual{\hhdivz}, K^{1/2} H_h^1 + L^2_h \right) = \mathcal{L}(\dual{\VV_h}, Q_h).
}
Since $\lif \grad_h$ is the identity on $Q_h$, we get from \eqref{eq:lifting-mapping-property} that for every $q \in Q_h$,
\eqns{
	\algnd{
		\norm{q}_{Q_h} &\leq C \norm{\grad_h q}_{\dual{\VV_h}} \\
		&= C \sup_{\vv \in \VV_h}\frac{\inner{q}{\div \vv}}{\norm{\vv}_{\VV_h}},
	}
}
where $C = \max(1, \beta^{-1})$ and is thus independent of both $h$ and $K$.  
\end{proof}

\TT{%
\begin{corollary}[A $K$ and $h$-robust preconditioner for the discrete mixed Darcy problem]\label{prop:disc-precon}
Define $\mathcal{A}_h$ as the coefficient matrix characterizing the left-hand 
side of the discrete problem \eqref{eqn:mixeddarcy:1:d}--%
\eqref{eqn:mixeddarcy:2:d}.  Then the preconditioner defined by %
\begin{equation}\label{eqn:disc_precond}
	\mathcal{B}_h =
	\begin{bmatrix}
		\left( K^{-1}I_h - \grad_h \div \right)^{-1} & 0 \\
		0 & I_h\TT{^{-1}} + (-K \div \grad_h)^{-1}
	\end{bmatrix}.
\end{equation}
provides a robust preconditioner for \eqref{eqn:mixeddarcy:1:d}--%
\eqref{eqn:mixeddarcy:2:d}.  That is, the condition number of 
$\mathcal{B}_h\mathcal{A}_h$ is bounded uniformly with respect to both $K$ 
and $h$.
\end{corollary}
}

\begin{proof}
\TT{%
The proof again follows the ideas of \cite{mardal2011preconditioning} and is 
nearly identical to the proof of Corollary \ref{prop:conts-precon}.  We 
therefore outline only the main ideas here.  The arguments in the proof 
of Proposition \ref{prop:uniform-inf-sup} imply that $\mathcal{A}_h$ is a 
homeomorphism from $\VV_h \times Q_h$ to its dual.  }%
%
%
\TT{%
Moreover, the norms on $\mathcal{A}_h$ and its inverse are bounded 
independently of $K$ and $h$.  Arguments analogous to those defining 
$\mathcal{B}$, in Section \ref{sec:uniform_conts}, lead directly to the 
framework \cite{mardal2011preconditioning} preconditioner given by the 
operator $B_h: \dual{\VV_h} \times \dual{Q_h} \to \VV_h \times Q_h$  
defined by \eqref{eqn:disc_precond}.
}
\end{proof}


\TT{%
\begin{remark}
We remark that a small liberty has been taken for the notation of $I_h$ 
in \eqref{eqn:disc_precond}.  Specifically, $I_h$ in the top left block of 
\eqref{eqn:disc_precond} signifies the identity $I_h:\VV_h\rightarrow \VV_h$ 
while the use of the same symbol in the bottom right block signifies the 
identity on $Q_h$; recall that the discrete
gradient, $\grad_h$, is defined by \eqref{eqn:gradh}.
\end{remark}
}

\section{Numerical Experiments}\label{sec:numerical-experiments}
Let $\Omega$ be a triangulation of the unit square such that the unit square is 
first divided in $N\times N$ squares of length $h=1/N$. Each square is then 
divided into two triangles. Below we will consider the case of homogeneous 
Dirichlet conditions for the flux; fluxes will be discretized by zero'th order 
RT elements.  In addition, we approximate the pressure in the space of 
piecewise constants and compute the eigenvalues of the preconditioned system. 
\TT{In the sections that follow, three different paradigms 
are presented: in Section \ref{subsec:numerical-experiments:spatially-constant} 
a spatially constant permeability, pursuant to the theory of 
Section \ref{sec:uniform_conts}, is demonstrated; in Section 
\ref{subsec:numerical-experiments:anisotropic-jump} we consider a computational 
experiment extending a jump in the (scalar) coefficient, $K$, and an 
anisotropic conductivity tensor; finally, in Section 
\ref{subsec:numerical-experiments:coupled-darcy-subsystem}, we compare the 
analogues of the preconditioning strategies analyzed here in the context 
of a coupled problem with Darcy subsystem.}

\subsection{\TT{A spatially constant conductivity}}
\label{subsec:numerical-experiments:spatially-constant}
Order-optimal multilevel methods for both $\hdiv$ and $H^1$ problems are 
well-known; thus, we consider preconditioners based simply on exact inversion.
It was shown in ~\cite{rusten1996interior} that the following local 
$-\Delta_h$ operator is spectrally equivalent to $-\div\nabla_h$ for discontinuous 
Lagrange elements of arbitrary order: 
\[
(-\Delta_h p, q) = \sum_{T\in \mathcal{T}_h} \int_T \nabla p \cdot \nabla q \, dx +  \sum_{E_i\in\mathcal{E_I}} \int_{E_i} \avg{h}^{-1}\jump{p}\jump{q} \, dS  +  \sum_{E_i\in\mathcal{E_D}} \int_{E_i} h^{-1} p q \, dS \quad\forall p, q\in Q_h.    
\]
Here, $\mathcal{T}_h$ is a triangulation of the domain $\Omega$; internal faces are signified by the set $\mathcal{E_I}$ whereas $\mathcal{E_D}$ denotes faces at the boundary associated with a pressure Dirichlet condition.
Furthermore $\avg{f}=\tfrac{1}{2}(f|_{T^{+}}+f|_{T^{-}})$ and
$\jump{f}=f|_{T^{+}}-f|_{T^{-}}$ are respectively the average and jump value of $f$
from the two elements $T^{\pm}$ that share the internal facet.  In Table \ref{tab:num_exper}  
we compare the discrete versions of the following two preconditioners 
\[
	\mathcal{B}_1 = \left[ \begin{array}{cc} (\frac{1}{K}(I - \nabla \div))^{-1} & 0 \\ 
           0 & (K I)^{-1}  \end{array} \right], 
\]
and  
\[
\mathcal{B}_2 = \left[ \begin{array}{cc} (\frac{1}{K}I - \nabla \div)^{-1} & 0 \\ 
           0 & \TT{I^{-1}} + (-K\Delta)^{-1} \end{array} \right]. 
\]

Table \ref{tab:num_exper} shows that both $\mathcal{B}_1$ and $\mathcal{B}_2$  
yield robust results for any $K \in (0,1)$, but that  $\mathcal{B}_1$ is 
usually somewhat better. This is somewhat surprising since the opposite is the 
case in simple tests with random matrices.  Numerical experiments (not reported 
here) confirm that the robust behavior applies also to Neumann conditions 
and BDM elements.  

\begin{table}
    \footnotesize{
      \begin{minipage}{0.49\textwidth}
          \begin{center}
            \begin{tabular}{c|c||cccc} \hline
              & \backslashbox{$K$}{$h$} &   $2^{-2}$  & $2^{-3}$ & $2^{-4}$ & $2^{-5}$   \\
              \hline  
    \multirow{5}{*}{$\mathcal{B}_1$} 
&        1 & 1.1 & 1.1 & 1.1 & 1.1\\
& $1^{-2}$ & 1.1 & 1.1 & 1.1 & 1.1\\
& $1^{-4}$ & 1.1 & 1.1 & 1.1 & 1.1\\
& $1^{-6}$ & 1.1 & 1.1 & 1.1 & 1.1\\
& $1^{-8}$ & 1.1 & 1.1 & 1.1 & 1.1\\
\hline
\multirow{5}{*}{$\mathcal{B}_2$}     
  & 1        & 1.1(2.4) & 1.1(1.4) & 1.1(1.4)& 1.1(1.4)\\
  & $10^{-2}$ & 2.3(2.3) & 2.3(2.3) & 2.3(2.3) & 2.3(2.3)\\
  & $10^{-4}$ & 2.9(2.6) & 3.1(2.6) & 3.1(2.6) & 3.0(2.6)\\
  & $10^{-6}$ & 3.0(2.6) & 3.3(2.6) & 3.5(2.6) & 3.5(2.6)\\
  & $10^{-8}$ & 3.0(2.6) & 3.3(2.6) & 3.5(2.6) & 3.6(2.6)\\
\hline
\end{tabular}
\end{center}
      \end{minipage}
    }
          \footnotesize{
            \begin{minipage}{0.49\textwidth}
              \begin{center}
                \begin{tabular}{c|c||cccc}
                  \hline
                  & \backslashbox{$K_{0}$}{$h$} &   $2^{-2}$  & $2^{-3}$ & $2^{-4}$ & $2^{-5}$   \\
                  \hline
    \multirow{5}{*}{$\mathcal{B}_1$} 
&        1 & 1.1 & 1.1 & 1.1 & 1.1\\
& $1^{-2}$ & 1.1 & 1.1 & 1.1 & 1.1\\
& $1^{-4}$ & 1.1 & 1.1 & 1.1 & 1.1\\
& $1^{-6}$ & 1.1 & 1.1 & 1.1 & 1.1\\
& $1^{-8}$ & 1.1 & 1.1 & 1.1 & 1.1\\
\hline
\multirow{5}{*}{$\mathcal{B}_2$} 
& 1        & 1.1(1.4) & 1.1(1.4) & 1.1(1.4) & 1.1(1.4)\\
& $10^{-2}$ & 2.5(2.1) & 2.4(2.1) & 2.4(2.1) & 2.4(2.1) \\
& $10^{-4}$ & 3.4(2.6) & 3.2(2.6) & 3.1(2.6) & 3.1(2.6) \\
& $10^{-6}$ & 3.4(2.6) & 3.5(2.6) & 3.4(2.6) & 3.4(2.6) \\
& $10^{-8}$ & 3.4(2.6) & 3.5(2.6) & 3.5(2.6) & 3.5(2.6) \\
\hline
                \end{tabular}
  \end{center}
            \end{minipage}
            }
          \caption{Condition numbers of the operators $\mathcal{B}_1 \mathcal{A}$ and
            $\mathcal{B}_2 \mathcal{A}$.
            (left) $K$ is constant. (right) $K$ is discontinuous given by \eqref{eq:jump_K}.
            Condition number of using the exact Schur complement (i.e. not $-\Delta_h$)
            is shown in the braces.
          }
          \label{tab:num_exper} 
\end{table}

\subsection{\TT{An anisotropic conductivity tensor with jump}}
\label{subsec:numerical-experiments:anisotropic-jump}
We now consider the case where a jump is present in the (scalar) coefficient, $K$, 
and the 
\TT{conductivity tensor} is given by an anisotropic matrix. These cases are not covered by the theoretical analysis and it is as such interesting
to compare the alternative preconditioners. To this end let $\Omega$ be a unit square and
\begin{equation}\label{eq:jump_K}
K=\begin{cases}
1 & x < \tfrac{1}{2}\\
K_0 & \mbox{otherwise}
\end{cases},
\quad\quad
    {\bf K}= {\bf R}_{\theta}{\bf \mbox{diag}}(1, K_0){\bf R}^{T}_{\theta},
\end{equation}
where $K_0 \in\left(0, 1\right]$ and ${\bf R}_{\theta}$ is a matrix of rotation
  by angle $\theta$.
  Homogeneous Dirichlet conditions for the normal flux shall be imposed on
  the left and right edges. As in the previous
  experiments we consider a uniform triangulation of the domain and 
  finite element discretization in terms of zero'th order RT and piecewise
  constant elements.

  Varying the magnitude of the $K_0$ Table \ref{tab:num_exper} shows
  that both preconditioners are practically unaffected by the presence
  of the considered discontinuity. We remark that the pressure preconditioner in
  the discrete $\mathcal{B}_2$ operator is discretized as follows
  \[
  \begin{split}
	  (-\mbox{div}(K \nabla_h p), q) =
    \sum_{E\in\mathcal{E}_{\mathcal{I}}} \int_{E_i}\frac{1}{\avg{h}}\frac{1}{\avg{K^{-1}}}\jump{p}\jump{q}dS+
  \sum_{E\in\mathcal{E}_{\mathcal{D}}} \int_{E_i}\frac{K}{h} p q dS\quad\forall p, q\in Q_h.
  \end{split}
  \]
  Note in particular that $K$ is averaged using the harmonic mean. 
  Table \ref{tab:tensor} shows the condition 
  numbers of the preconditioned %
  problems with matrix valued permeability. Here the preconditioners $\mathcal{B}_1$ and
  $\mathcal{B}_2$ are generalized as
  \[
  \mathcal{B}_1=\begin{bmatrix}
  ({\bf K}^{-1}I - \nabla\div{\bf K}^{-1})^{-1} & 0\\
    0 & (kI)^{-1}
  \end{bmatrix}\mbox{ and }
  \mathcal{B}_2=\begin{bmatrix}
  ({\bf K}^{-1}I - \nabla\div)^{-1} & 0\\
    0 & I^{-1} + (-\div(k\nabla))^{-1}
    \end{bmatrix}
  \]
  with $k=1/\sum_i\lambda_i({\bf K})^{-1}$ and $\lambda_i({\bf K})$ the $i$-th
  eigenvalue of the permeability matrix. 

\miro{%
\begin{remark}
 We note that the leading
 block of $\mathcal{B}_1$ is implemented as a solver for the variational 
 problem.  Namely, we find $\uu\in \hdivz$ such that
 \[
    (\Div\left({\bf K}^{-1}\uu\right),\Div \vv) + ({\bf K}^{-1}\uu, \vv) = (\mathbf{f}, \vv), \quad\forall\vv\in\hdivz.
    \]
\end{remark}
}
  
  The preconditioner $\mathcal{B}_1$ is robust
  with respect to the conditioning of ${\bf K}$; in particular the simple scaling of
  the pressure mass matrix is sufficient. This stands in contrast to the case of 
  preconditioner $\mathcal{B}_2$; here, the same scaling yields $h$ independent 
  condition numbers only for relatively large $K_0$.
  In order to recover the mesh and parameter robustness attributes for $\mathcal{B}_2$ 
  we approximate the operator $-\div({\bf K}\nabla_h)$ as the exact Schur
  complement of $\mathcal{A}_h$. That is, an exact inverse of ${\bf K}^{-1}I_h$ on $V_h$
  is used in the construction of the preconditioner; c.f.~Table \ref{tab:tensor}. 
  We remark that \cite{powell2003optimal} discusses a construction based on a 
  diagonalized mass matrix. In particular, for ${\bf K}$ diagonal, a Jacobi 
  preconditioner is shown to yield bounds independent of ${\bf K}$.

  \begin{table}[h]
    \begin{minipage}{0.49\textwidth}
      \centering
      \footnotesize{
              \begin{center}
          \begin{tabular}{c|c|c||cccc}
            \hline
 &           $\theta$ & \backslashbox{$K_0$}{$h$} & $2^{-3}$ & $2^{-4}$ & $2^{-5}$ & $2^{-6}$\\
            \hline
\multirow{10}{*}{$\mathcal{B}_1$}           
&\multirow{5}{*}{$0$}	&1	   &1.4	&1.4	&1.4	&1.4\\     
&                        &$10^{-2}$  &1.4	&1.4	&1.4	&1.4\\	
&                        &$10^{-4}$  &1.4	&1.4	&1.4	&1.4\\	
&                        &$10^{-6}$  &1.4	&1.4	&1.4	&1.4\\	
&                        &$10^{-8}$  &1.4	&1.4	&1.4	&1.4\\
\cline{2-7}
&\multirow{5}{*}{$\pi/4$} &1	   &1.4	&1.4	&1.4	&1.4\\     
&                         &$10^{-2}$ &1.4	&1.4	&1.4	&1.4\\	
&                         &$10^{-4}$ &1.4	&1.4	&1.4	&1.4\\
&                         &$10^{-6}$ &1.4	&1.4	&1.4	&1.4\\
&                         &$10^{-8}$ &1.4	&1.4	&1.4	&1.4\\
\hline
          \end{tabular}
              \end{center}
              }
  \end{minipage}
    \begin{minipage}{0.49\textwidth}
      \centering
          \footnotesize{
                \begin{center}
          \begin{tabular}{c|c|c||ccc}
            \hline
          &  {$\theta$} & \backslashbox{$K_0$}{$h$} & $2^{-3}$ & $2^{-4}$ & $2^{-5}$\\
            \hline
\multirow{10}{*}{$\mathcal{B}_2$}            
&\multirow{5}{*}{$0$}     &1	    &1.1 (1.1)	&1.1 (1.1)	&1.1 (1.1) \\     
&                         &$10^{-2}$ &3.4 (2.4) &3.4 (2.4)      &3.4 (2.4) \\	
&                         &$10^{-4}$ &13.3(2.6)	&21.0(2.6)	&27.7(2.6) \\
&                         &$10^{-6}$ &14.4(2.6)	&27.2(2.6)	&52.8(2.6) \\
&                         &$10^{-8}$ &14.5(2.6)	&27.3(2.6)	&53.5(2.6) \\
\cline{2-6}                                                                          
&\multirow{5}{*}{$\pi/4$}	&1	   &1.1  (1.1)       &1.1  (1.1)    &1.1 (1.1) \\     
&                        &$10^{-2}$ &3.9  (1.9)       &4.2  (1.7)    &4.3 (1.5) \\	
&                        &$10^{-4}$ &11.5 (2.6)       &19.3 (2.6)    &27.9(2.6) \\	
&                        &$10^{-6}$ &12.5 (2.6)       &25.0 (2.6)    &50.2(2.6) \\	
&                        &$10^{-8}$ &12.5 (2.6)       &25.1 (2.6)    &50.9(2.6) \\
\hline
          \end{tabular}
                \end{center}
                }
  \end{minipage}
    \caption{Condition numbers of the operators $\mathcal{B}_1\mathcal{A}$ and
      $\mathcal{B}_2\mathcal{A}$ with permeability matrix ${\bf K}$ in \eqref{eq:jump_K}.
      Numbers in parentheses denote the condition number of $\mathcal{B}_2\mathcal{A}$ 
      when using the exact Schur complement to approximation $-\div({\bf K}\nabla)$.
  }
  \label{tab:tensor}
\end{table}

\subsection{\TT{A coupled problem with a Darcy subsystem}}
\label{subsec:numerical-experiments:coupled-darcy-subsystem}
\kent{To illustrate the potential benefit of having the possibility to choose 
different scales of the pressure, we consider a simplified Biot problem.}  
\TT{We remark that fully-parameter-robust preconditioners for Biot's problem 
are established, in detail, in \cite{hong2017parameter}; the discussion 
here is intended only to illustrate the difficulties encountered in a single 
parameter regime}.  \TT{With this caveat in mind: consider a simplified, 
three-field formulation of Biot's problem that has been differentiated in time 
via an implicit Euler scheme; the resulting system has the form:}
\TT{%
  \begin{equation}\label{eq:simple_biot}
\begin{aligned}
-\div\sigma(\uu) + \grad p &= \mathbf{f}_{\uu}, &\mbox{ in }\Omega,\\
\frac{1}{K} \vv + \grad p &= \mathbf{f}_{\vv}, &\mbox{ in }\Omega,\\
- \div\uu - \div\vv &= f_p, &\mbox{ in }\Omega,%
\end{aligned}%
\end{equation}
}
\TT{where $\sigma(\uu) = 2\epsilon(\uu) = \left(\grad \uu + \grad \uu^T\right)$ and 
$K\in (0,1]$ is assumed fixed, and spatially constant, but otherwise 
arbitrary.} 
\miro{Here the system \eqref{eq:simple_biot} is completed by the 
    Dirichlet boundary conditions on the displacement $\uu$, i.e. $\uu=\uu_0$
on $\Gamma_D\subset\partial\Omega$, and Neumann condition $\sigma\cdot\mathbf{n}=\mathbf{h}_0$ on
$\Gamma_N=\partial\Omega\setminus\Gamma_D$. For the Darcy subproblem we then set $\vv\cdot \mathbf{n} = \vv_0$ on
$\Gamma_D$ and $p=p_0$ on $\Gamma_N$.}
\kent{It is here clear that the pressure provides the coupling between the elastic deformation and the porous
media flow and that Brezzi theory of the two sub-systems of Stokes and Darcy type
enables stability of the coupled system~\cite{hong2017parameter,rodrigo2018new}.  
The Stokes problem requires
that the displacement is bounded in $\mathbf{H}^1$,  while  the pressure is bounded in $L^2$. 
For the Darcy problem, a flux in $\hdivz$ can be readily combined with a pressure in $L^2$. However,  
when $K\rightarrow 0$, the norms must be scaled. An in particular, as 
proposed in~\cite{powell2003optimal,vassilevski2013block,vassilevski2014mixed}  
a natural norm for the pressure is $\sqrt K L^2$.
However, the Stokes-type coupling
between the displacement and pressure suggests the pressure should not be scaled with $K$. 
On the other hand, the pressure norm $L^2+ K^{1/2}$, used in our paper, can be bounded by $L^2$. 
Below, we illustrate some of the problems that may occur in a preconditioning setting 
with  wrongly scaled pressures and demonstrate that the scaling proposed in this paper
can indeed be combined into a robust preconditioner for a simplified Biot problem. }

\miro{As preconditioners for \eqref{eq:simple_biot} we shall consider operators
\begin{equation}\label{eq:B_biot}
  \mathcal{B}_1=\begin{pmatrix}
  -\Div\epsilon & &\\
  & K^{-1}I -K^{-1}\nabla\Div & \\
  & & I\\
  \end{pmatrix}^{-1}\mbox{ and }
  \mathcal{B}_2=\begin{pmatrix}
  -\Div\epsilon & &\\
  & K^{-1}I -\nabla\Div & \\
  & & I\\
  \end{pmatrix}^{-1},
\end{equation}
where, importantly, $\mathcal{B}_2$ avoids rescaling of the $\nabla\Div$ block,
cf. \eqref{their_prec:a} and the related discussion of the issue in the context of
multiphysics problems. In this respect, the
operator is similar to the Darcy preconditioner \eqref{eq:continuous-preconditioner}.
Our aim is then to demonstrate that such scaling may be not be desirable for
$K$-robustness. To this end we perform a numerical experiment where $\Omega=(0, 1)^2$ with $\Gamma_D$
the union of left, right and bottom edge of $\Omega$. Furthermore a
uniform triangulation of the domain is considered with spaces of the discrete displacement,
flux $\vv$ and the pressure constructed respectively by continuous piecewise quadratic
Lagrange elements, zero'th order RT elements and piecewise constant
elements. Indeed, Table \ref{tab:biot_cond} shows that only the preconditioner $\mathcal{B}_2$
leads to condition numbers bounded in $K$.}
\TT{One may cautiously observe that the pressure block in \eqref{their_prec:a} 
includes a $K$-scaling.  It is a reasonable enquiry, then, to ask whether 
including this scaling can rectify the poor performance observed 
for the $\mathcal{B}_1$ preconditioner.  To test this observation} %
%
\miro{let us also consider preconditioning \eqref{eq:simple_biot} with operators
\begin{equation}\label{eq:B_biot_K}
  \mathcal{B}^K_1=\begin{pmatrix}
  -\Div\epsilon & &\\
  & K^{-1}I -K^{-1}\nabla\Div & \\
  & & KI\\
  \end{pmatrix}^{-1}\mbox{ and }
  \mathcal{B}^K_2=\begin{pmatrix}
  -\Div\epsilon & &\\
  & K^{-1}I -\nabla\Div & \\
  & & KI\\
  \end{pmatrix}^{-1}.
\end{equation}
}
%
\TT{Condition numbers corresponding to these operators are shown in} %
\miro{Table \ref{tab:biot_cond}.
It is evident that neither of the preconditioners of 
\eqref{eq:B_biot_K} present any improvement compared to
their unscaled variants. We conclude, comparing $\mathcal{B}^K_2$ and $\mathcal{B}_2$,
that rescaling the pressure is not necessary for robustness with respect to $K$.
}

\begin{table}[h!]
  \centering
  \footnotesize{
    \begin{minipage}{0.4\textwidth}
      \centering
    \begin{tabular}{c|c||lll}
\hline      
 &{\backslashbox{$K$}{$h$}} & $2^{-2}$ & $2^{-3}$ & $2^{-4}$\\
\hline
\multirow{7}{*}{$\mathcal{B}_1$} & $1$      &3.5	        &4.0	&4.0  \\
&$10^{-1}$ &14.3	&18.2	&19.0 \\
&$10^{-2}$ &91	        &118	&139 \\
&$10^{-3}$ &396	        &731	&$1.0\times10^{3}$  \\
&$10^{-4}$ &739	        &$2.5\times10^{3}$	&$5.2\times10^{3}$  \\
&$10^{-6}$ &817	        &$3.6\times10^{3}$	&$14\times10^{3}$ \\
&$10^{-8}$ &818	        &$3.6\times10^{3}$	&$15\times10^{3}$ \\
\hline
\multirow{7}{*}{$\mathcal{B}^K_1$} & $1$      &3.5	&4.0	&4.0\\
&$10^{-1}$ &29.8	&39.6	        &41.8\\
&$10^{-2}$ &559	        &744	        &893\\
&$10^{-3}$ &7$\times 10^{3}$	&$14\times10^{3}$	        &$20\times10^{3}$\\
&$10^{-4}$ &4$\times 10^{4}$	&$15\times10^{4}$	        &$32\times10^{4}$\\
&$10^{-6}$ &5$\times 10^{5}$	&$22\times10^{5}$	        &$89\times10^{5}$\\
&$10^{-8}$ &5$\times 10^{6}$	&$22\times10^{6}$	        &$90\times10^{6}$\\
\hline
    \end{tabular}
    \end{minipage}
    \begin{minipage}{0.58\textwidth}
      \centering
    \begin{tabular}{c|c||lllll}
\hline
& {\backslashbox{$K$}{$h$}} & $2^{-2}$ & $2^{-3}$ & $2^{-4}$ & $2^{-5}$ & $2^{-6}$\\
\hline
\multirow{7}{*}{$\mathcal{B}_2$} & $1$      &3.5	        &3.8	&4.0	&4.0	&4.0    \\
&$10^{-1}$ &5.7         &5.9	&6.0	&6.0	&6.0    \\
&$10^{-2}$ &11.0	&11.4	&11.6	&11.8	&11.8   \\
&$10^{-3}$ &12.2	&12.9	&13.4	&13.6	&13.8   \\
&$10^{-4}$ &12.1	&12.9	&13.3	&13.7	&14.1   \\
&$10^{-6}$ &12.1	&12.8	&13.1	&13.3	&13.3   \\
&$10^{-8}$ &12.1	&12.8	&13.1	&13.3	&13.3   \\
\hline
\multirow{7}{*}{$\mathcal{B}^K_2$}
&$1$       &3.5	&3.80	&4.0	&4.0	&4.0\\                                   
&$10^{-1}$ &7.8	&8.7	&9.2	&9.5	&9.6\\                                   
&$10^{-2}$ &19.8 &23.0	&25.6	&27.2	&28.0\\                                  
&$10^{-3}$ &40	&53	&66	&74	&81\\                                  
&$10^{-4}$ &95	&114	&141	&182	&132\\                                 
&$10^{-6}$ &895	&960	&995	&1029	&1114\\                        
&$10^{-8}$ &8947 &9575	&9854	&9946	&9986\\
\hline
    \end{tabular}
    \end{minipage}
    
  }
  \caption{Condition numbers of \eqref{eq:simple_biot} with preconditioners
    \eqref{eq:B_biot} and \eqref{eq:B_biot_K}.
    For $\mathcal{B}_1$ and $\mathcal{B}^K_1$ only the results from direct solver (computing
    the entire spectrum) are included as iterative solver did not converge within
    1500 iterations.
  }
  \label{tab:biot_cond}
\end{table}

\section{Conclusion}\label{sec:conclusion}
In this paper we have introduced a new preconditioner for the Darcy problem and shown that for $K\in\reals+$ the preconditioner
is stable as $K\rightarrow 0$. 
The preconditioner is of the form
\[
\mathcal{B}_2 = \left[\begin{array}{cc}  \TT{\left(\alpha A + B^T B\right)^{-1}}  & 0 \\ 0 & \TT{I^{-1} +  (B(\alpha A)^{-1} B^T)^{-1}}    
\end{array} \right].
\]
We have compared it with a simpler approach where a scaling is introduced on the divergence term,
i.e. the preconditioner of the form: 
\[
\label{their_prec}
\mathcal{B}_1 = \left[\begin{array}{cc} \TT{(\alpha A + \alpha B^T B)^{-1}}  & 0 \\ 0 & \TT{(\frac{1}{\alpha} I)^{-1}} \end{array} \right].
\]
The preconditioners were tested with both  constant and small  permeabilities, i.e. from $10^{-8}$ to 1 and 
permeabilities that contain jumps and anisotropy of similar sizes.  Both preconditioners work well;
though $\mathcal{B}_1$ appears to result in a smaller condition number. 
$\mathcal{B}_1$ has the advantage of performing better, directly, in the case of matrix-valued $K$.
However, $\mathcal{B}_2$ may be of interest in multi-physics codes where scaling of the $B^TB$ component is not 
desired; in this case, significantly improved condition numbers have been observed when using the 
exact Schur complement to approximate $-\div(K\nabla)$ in the presence of matrix-valued $K$.  That is, in
the presence of anisotropy, a crucial component is a to find a proper local operator representing  
$B(\alpha A)^{-1} B^T$.

\section{Acknowledgements}
The authors gratefully acknowledge helpful discussions with Ragnar Winther\footnote{Dept.~of Mathematics, University of Oslo. Oslo, Norway.}
in addition to productive exchanges with Marie E.~Rognes\footnote{Dept.~of Numerical Anal.~and Sci.~Comput.~, Simula Research Laboratory. Fornebu, Norway} and Walter Zulehner\footnote{ Johannes Kepler University Linz, Institute of Computational Mathematics Linz, Austria}.  The work of Miroslav Kuchta was 
supported by the Research Council of Norway (NFR) grant 280709.  The work of Travis Thompson was supported by the Research Council of Norway 
under the FRINATEK Young Research Talents Programme through project \#250731/F20 (Waterscape).

\bibliographystyle{abbrv}
\bibliography{note}
\end{document}